\documentclass[12pt]{amsart}
\usepackage[mathscr]{eucal}
\usepackage{amssymb,amsmath}
\usepackage{amsfonts}
\usepackage{a4}

\usepackage{enumerate}

\theoremstyle{plain}
\newtheorem{thm}{Theorem}[section]
\newtheorem{prp}[thm]{Proposition}
\newtheorem{cor}[thm]{Corollary}
\newtheorem{lem}[thm]{Lemma}

\newtheorem{ex}[thm]{Example}

\theoremstyle{remark}
\newtheorem{prob}[thm]{Problem}
\newtheorem{rem}[thm]{Remark}

\numberwithin{equation}{section}

\DeclareMathOperator{\supp}{supp}
\DeclareMathOperator{\clo}{Cl}

\newcommand{\co}{c_0(\Gamma)}

\newcommand{\phe}{\varphi}

\newcommand{\ec}{\mathcal{EC}}
\newcommand{\ecn}{\mathcal{EC}_n}
\newcommand{\ecf}{\mathcal{EC}_{<\omega}}
\newcommand{\eco}{\mathcal{EC}_{\omega c}}

\begin{document}

\title{On two problems concerning Eberlein compacta}

\author{Witold Marciszewski}
\address{Institute of Mathematics\\
University of Warsaw\\ Banacha 2\newline 02--097 Warszawa\\
Poland}
\email{wmarcisz@mimuw.edu.pl}

\thanks{The author was partially supported by the NCN (National Science Centre, Poland) research grant no.\ 2020/37/B/ST1/02613.}

\subjclass[2020]{Primary 46A50, 54D30, 54G12}
\keywords{Eberlein compact, Corson compact, scattered space, $C(K)$}

\begin{abstract} 
We discuss two problems concerning the class Eberlein compacta, i.e., weakly compact subspaces of Banach spaces. The first one deals with preservation of some classes of scattered Eberlein compacta under continuous images. The second one concerns the known problem of the existence of nonmetrizable compact spaces without nonmetrizable zero-dimen\-sional closed subspaces. We show that the existence of such Eberlein compacta is consistent with \textsf{ZFC}. We also show that it is consistent with \textsf{ZFC} that each Eberlein compact space of weight $> \omega_1$ contains a nonmetrizable closed  zero-dimensional subspace.
\end{abstract}

\maketitle

\section{Introduction}\label{intro}
This paper is concerned with Eberlein compacta, i.e., compact spaces which can be embedded into a Banach space equipped with the weak topology. This class of spaces has been intensively studied for its interesting topological properties and various connections to functional analysis; we refer the reader
to a survey article by Negrepontis \cite{Ne}. 

It is well known that this class of compacta and its several subclasses, like classes of uniform Eberlein compacta, scattered Eberlein compacta, scattered Eberlein compacta of height $\le n$ (see Section 2 for definitions), are preserved by continuous images. We will discuss the problem of preservation under continuous images for some classes of scattered Eberlein compacta  $K$ closely related to the properties of the Banach space $C(K)$ of continuous real-valued functions on $K$. 

For a set $X$ and $n\in\omega$, by $\sigma_{n}(X)$ we denote the subspace of the product $2^X$
consisting of all functions with supports of cardinality $\le n$. Given $n\in\omega$, we will say that a compact space $K$ belongs to the class $\ecn$ if $K$ can be embedded in the space
$\sigma_{n}(X)$ for some set $X$. The class $\ecn$ is a proper subclass of the class of
scattered Eberlein compacta  of height $n+1$. We will denote the union $\bigcup_{n\in\omega}\ecn$ by $\ecf$.

In \cite{Ma} it was proved that, for a compact space $K$, the Banach space $C(K)$ is isomorphic  to $c_0(\Gamma)$, for some set $\Gamma$, if and only if, the space $K$ belongs to $\ecf$, see Theorem \ref{old_charact}. From this characterization we derive that the class  $\ecf$ is preserved by continuous images, see Corollary \ref{cor_images}. However, we show that this does not hold true for the classes $\ecn$. We give an example of a continuous image $L$ of a space $K\in \ec_2$ such that $L$ does not belong to $\ec_2$ (Example \ref{ex_cont_image}). We also prove that each continuous image of a space from $ \ec_2$ belongs to $\ec_3$ (Theorem \ref{k(2)}).
In general case, we show that, for each $n\in\omega$, there exists $k(n)\in\omega$ such that any continuous image of a space $K\in \ecn$ belongs to $\ec_{k(n)}$ (Theorem \ref{k(n)}). These results are can be found in Section 3.

The last section of the paper is devoted to the known problem of the existence of nonmetrizable compact spaces without nonmetrizable zero-dimen\-sional closed subspaces. Several such spaces were obtained using 
some additional set-theoretic assumptions. Recently,
P.\ Koszmider \cite{Ko} constructed the first such example in \textsf{ZFC}. We investigate this problem for the class of Eberlein compact spaces.
We construct such Eberlein compacta, assuming the existence of a Luzin set, see Corollary \ref{ex_no_zerodim_sub}. We also show that it is consistent with \textsf{ZFC} that each Eberlein compact space of weight $\ge \omega_2$ contains a nonmetrizable closed  zero-dimensional subspace (Corollary \ref{cor_corson}).

\section{Terminology and some auxiliary results}\label{notation}

\subsection{Notation}

All topological spaces under consideration are assumed to be Tikhonov.

For a set $X$ and $n\in\omega$, we use the standard notation
$[X]^n=\{A\subseteq X: |A|=n\}$, $[X]^{\le n}=\bigcup\{[X]^k: k\le
n\}$ and $[X]^{< \omega}=\bigcup\{[X]^k: k< \omega\}$. 

We say that a family $\mathcal{U}$ of sets has \emph{order $\le n$} if
every subfamily
$\mathcal{V}\subset\mathcal{U}$ of cardinality $n+2$ has an empty intersection (in other terminology, the family $\mathcal{U}$ is
\emph{point-$(n+1)$}). The family $\mathcal{U}$ has \emph{finite order} if it has order $\le n$ for some $n\in\omega$.

The family $\mathcal{U}$ of subsets of a space $X$
is $T_0$-separating if, for every pair of distinct points $x,y$ of
$X$, there is $U\in\mathcal{U}$ containing exactly one of the
points $x,y$.

For a locally compact space $X$, $\alpha(X)$ denotes the one point compactification of $X$. We denote the point at infinity of this compactification by $\infty_X$. 

\subsection{Function spaces}

Given a compact space $K$, by $C(K)$ we denote the Banach space of
continuous real-valued functions on $K$, equipped with the standard
supremum norm. 

\subsection{Scattered spaces}

A space $X$ is \textit{scattered} if no nonempty subset $A\subseteq X$ is dense-in-itself.

For a scattered space $K$, by
\textit{Cantor-Ben\-dixson height} $ht(X)$ of $K$ we mean the minimal ordinal
$\alpha$ such that the Cantor-Bendixson derivative $K^{(\alpha)}$
of the space $K$ is empty. The Cantor-Ben\-dixson height of a
compact scattered space is always a nonlimit ordinal.

A surjective map $f:X\to Y$ between topological spaces is said to be \textit{irreducible} if no proper closed subset of $X$ maps onto $Y$.
If $X$ is compact, by Kuratowski-Zorn Lemma,
for any surjective map $f:X\to Y$, there is a closed subset $C\subseteq X$ such that the restriction $f\upharpoonright C$ is irreducible.

The following facts concerning continuous maps of scattered compact spaces are well known, cf.\ the proof of Proposition 8.5.3 and Exercise 8.5.10(C) in \cite{Se}.

\begin{prp}\label{scatt_Fact1}
Let $K$ be a scattered compact space and let $\phe: K\to L$ be a continuous surjection. Then, for each ordinal $\alpha$, we have $L^{(\alpha)} \subseteq \phe(K^{(\alpha)})$. In particular, $ht(L) \le ht(K)$.
\end{prp}

\begin{prp}\label{scatt_Fact2}
Let $K$ be a scattered compact space and let $\phe: K\to L$ be a continuous irreducible surjection. Then $L' = \phe(K')$ and $\phe\upharpoonright (K\setminus K')$ is a bijection onto $L\setminus L'$.
\end{prp}

\subsection{Eberlein and Corson compact spaces}

A space $K$ is an \emph{Eberlein} compact space
if $K$ is homeomorphic to a weakly compact subset of a Banach
space. 
Equivalently, a compact space $K$ is an Eberlein compactum if $K$ can 
be embedded in the following subspace of the product $\mathbb{R}^\Gamma$:
$$c_0(\Gamma)=\{x\in \mathbb{R}^\Gamma:\text{ for every $\varepsilon>0$ the set $\{\gamma: 
|x(\gamma)|>\varepsilon\}$ is finite}\},$$
for some set $\Gamma$, see \cite{Ne}.

If $K$ is homeomorphic to a weakly compact subset of a
Hilbert space, then we say that $K$ is a \emph{uniform Eberlein}
compact space. All metrizable compacta are uniform Eberlein.

A compact space $K$ is {\em Corson compact} if, for some set $\Gamma$, 
$K$ is homeomorphic to a subset of the $\Sigma$--product of real lines
$$\Sigma(\mathbb{R}^\Gamma)=\{x\in \mathbb{R}^\Gamma: |\{\gamma: x(\gamma)\neq 0\}|\le\omega\}.$$

Clearly, the class of Corson compact spaces contains all Eberlein compacta.

\subsection{Spaces $\sigma_{n}(X)$}\label{subsec_sigma_n}

Given a set $\Gamma$ and $n\in\omega$, by $\sigma_{n}(\Gamma)$ we denote the subspace of the product $2^\Gamma$
consisting of all characteristic functions of sets of cardinality $\le n$.  The space $\sigma_{n}(\Gamma)$ is uniform Eberlein and scattered of height $n+1$.

For $A\in [\Gamma]^{\le n}$, we denote
the standard clopen neighborhood $\{\chi_B\in
\sigma_{n}(\Gamma): A\subset B\}$ of $\chi_A$ in $\sigma_{n}(\Gamma)$ by $V_{A}$.

To simplify the notation we will say that a compact space $K$ belongs to the class $\ecn$ if $K$ can be embedded in the space
$\sigma_{n}(\Gamma)$ for some set $\Gamma$. We will denote the union $\bigcup_{n\in\omega}\ecn$ by $\ecf$. Trivially, if a compact space $K$ belongs to any of the above classes, then each closed subset of $K$ is also in the same class. One can also easily 
verify that the class $\ecf$ is preserved under taking finite products, cf.\ \cite[p.\ 148]{Av}.

\begin{prp}\label{charact_ecn}
For a compact space $K$ and $n\in\omega$, the following conditions are equivalent:
\begin{itemize}
\item[(i)] $K$ has a $T_0$-separating point-$n$ family of clopen
subsets;
\item[(ii)] $K$  belongs to the class $\ecn$.
\end{itemize}
\end{prp}

\begin{proof}
((i)$\Rightarrow$(ii)) Let $\mathcal{A}$ be a $T_0$-separating point-$n$ family of clopen
subsets of $K$. For $x\in K$, let $f_x: \mathcal{A}\to 2$  be a function defined by $f_x(A) = 1$ if $x\in A$, $0$ otherwise, for $A\in \mathcal{A}$. Clearly, the mapping $x\mapsto f_x$ is a required embedding.
\smallskip

\noindent
((ii)$\Rightarrow$(i)) Suppose that $K$ is a subspace of the space
$\sigma_{n}(\Gamma)$ for some set $\Gamma$. For $\gamma\in \Gamma$, let $U_\gamma = \{x\in K: x(\gamma) = 1\}$. One can easily verify that the family $\{U_\gamma: \gamma\in \Gamma\}$ is a $T_0$-separating point-$n$ family of clopen subsets of $K$.
\end{proof}

\begin{lem}\label{minimal_X}
Let $K$ be an infinite compact subset of $\sigma_{n}(\Gamma)$ for some set $\Gamma$ and $n\in\omega$. Then $K$ can can be embedded into $\sigma_{n}(\kappa)$, where $\kappa$ is the weight $w(K)$ of $K$.  
\end{lem}

\begin{proof}
Follows from the proof of Lemma \ref{charact_ecn} and the well known fact that, for an infinite compact space the cardinality of the family of clopen subsets of $K$ is bounded by $w(K)$.
\end{proof}

\begin{lem}\label{emb_n_into_n+1}
Let $\Gamma$ be an infinite set. Then for any $n,k\in \omega,\ k\ge 1$, the discrete union of $k$ copies of $\sigma_{n}(\Gamma)$ embeds into $\sigma_{n+1}(\Gamma)$.
\end{lem}

\begin{proof}
Let $X=\{x_0,x_1,\dots,x_{k-1}\}$ be a set disjoint with $\Gamma$. For $f\in \sigma_{n}(\Gamma)$ and $i< k$ let $f_i: \Gamma\cup X\to 2$ be defined by 
\begin{eqnarray*}
f_i(x) = 
\begin{cases}
f(x)& \text{if } x\in \Gamma\\
1& \text{if } x = x_i\\
0& \text{if } x = x_j,\ j<k,\ j\ne i
\end{cases}
\end{eqnarray*}

One can easily verify that, if we assign to a function $f$ from $i$-th copy of $\sigma_{n}(\Gamma)$, the function $f_i$, then we will obtain an embedding of the discrete union of $k$ copies of $\sigma_{n}(\Gamma)$ into $\sigma_{n+1}(\Gamma\cup X)$, a copy of $\sigma_{n+1}(\Gamma)$.
\end{proof}

\begin{thm}[Argyros and Godefroy] \label{thmAG}
Every Eberlein compactum $K$ of weight $<\omega_\omega$ and of
finite height belongs to the class $\ecf$.
\end{thm}

\begin{ex}[Bell and Marciszewski \cite{BM}]\label{exBM}
There exists an Eberlein compactum $K$ of weight $\omega_\omega$
and height $3$ which does not belong to $\ecf$.
\end{ex}

\subsection{Luzin sets and its variants}\label{subs_Luzin} 
Usually, a subset $L$ of real line $\mathbb{R}$ is called a \emph{Luzin set} if $X$ is uncountable and, for every meager subset $A$ of $\mathbb{R}$ the intersection $A\cap L$ is countable.
Let  $\kappa \le \lambda$ be uncountable cardinal numbers. We will say that a subset $L$ of a Polish space $X$ is a \emph{$(\lambda,\kappa)$-Luzin set} if $X$ has the cardinality $\lambda$, and, for every meager subset $A$ of $X$ the intersection $A\cap L$ has the cardinality less than $\kappa$. 
In this terminology, the existence of a Luzin set in $\mathbb{R}$ is equivalent to the existence of a $(\omega_1,\omega_1)$-Luzin set.

Since, for every Polish space $X$ without isolated points there is a Borel isomorphism $h:X\to \mathbb{R}$ such that $A\subseteq X$ is meager if and only if, $h(A)$ is meager in $\mathbb{R}$, it follows that the existence of a $(\lambda,\kappa)$-Luzin set in such $X$ is equivalent to the existence of a $(\lambda,\kappa)$-Luzin set in $\mathbb{R}$. 

It is known that, for each $n\ge 1$ the existence of a $ (\omega_n,\omega_1)$-Luzin set in $\mathbb{R}$ is consistent with \textsf{ZFC}, cf.\ \cite[Lemma 8.2.6]{BJ}.

\subsection{Cardinal numbers $\mathfrak{b}$ and $\mathrm{non}(\mathcal{M})$}\label{subs_b_nonM}

Recall that the preorder $\le^*$ on $\omega^\omega$ is defined by $f\le^* g$ if $f(n)\le g(n)$ for all but finitely $n\in\omega$. A subset $A$ of $\omega^\omega$ is called \emph{unbounded} if it is unbounded with respect to this preorder. In Section 4 we will use two cardinal numbers related with the structure of the real line
\begin{eqnarray*}
\mathfrak{b} &=& \min\{|A|: A \text{ is an unbounded subset of } \omega^\omega\}\\
\mathrm{non}(\mathcal{M})  &=& \min\{|B|: B \text{ is a nonmeager subset of } \mathbb{R}\}\,.
\end{eqnarray*}
It is well known that $\mathfrak{b}\le \mathrm{non}(\mathcal{M})$ (cf.\ \cite[Ch.\ 2]{BJ}), and, for each natural number $n\ge 1$, the statement $\mathfrak{b} = \omega_n$ is consistent with \textsf{ZFC}, (cf.\ \cite[Theorem 5.1]{vD}).

\subsection{Aleksandrov duplicate $AD(K)$ of a compact space $K$}\label{subs_AD}

Recall the construction of the  Aleksandrov duplicate $AD(K)$ of a compact space $K$. 

$AD(K) = K\times 2$, points $(x,1)$, for $x\in K$, are isolated in $AD(K)$ and basic neighborhoods of a point $(x,0)$ have the form $(U\times 2)\setminus \{(x,1)\}$, where $U$ is an open neighborhood of $x$ in $K$.

The following fact is well known (cf.\ \cite{KM}).

\begin{prp}\label{prp_AD_Eberlein}
The Aleksandrov duplicate $AD(K)$ of an (uniform) Eberlein compact space $K$ is (uniform) Eberlein compact.
\end{prp}

\begin{proof}
Without loss of generality we can assume that $K$ is a subspace of $c_0(\Gamma)$ ($\ell_2(\Gamma)$), equipped with the pointwise topology, for some set $\Gamma$. We will show that  $AD(K)$ can be embedded into the space $c_0(\Gamma\cup K)$ ($\ell_2(\Gamma\cup K)$). For $x\in K$ and $i=0,1$ define a function $f_{x,i}:\Gamma\cup K\to \mathbb{R}$ as follows:
$$f_{x,i}(t)=\left\{
    \begin{array}{llll}
        x(t) & \mbox{if } t\in \Gamma  \\
        0 & \mbox{if } t\in K,\; t\ne x \\
        1 & \mbox{if } t=x,\; i=1 \\
        0 & \mbox{if } t=x,\; i=0
    \end{array}
\right.$$
One can easily verify that the mapping $(x,i)\mapsto f_{x,i}$ gives the desired embedding.
\end{proof}

\section{On continuous images of compact subsets of spaces $\sigma_{n}(\Gamma)$}

The class $\ecf$ of compact subsets $K$ of spaces $\sigma_{n}(\Gamma), n\in\omega$ (cf.\ Subsection \ref{subsec_sigma_n}), turned out to be the class of those compacta $K$, for which the Banach space $C(K)$ is isomorphic to the Banach space $\co$ for some set $\Gamma$.

\begin{thm}[{\cite[Theorem 1.1]{Ma}}]
\label{old_charact} 
For a compact space $K$ the following
conditions are equivalent:
\begin{itemize}
\item[(i)] $K$ has a $T_0$-separating family of clopen
subsets and of finite order;
\item[(ii)] $K$  belongs to the class $\ecf$,
\item[(iii)] $C(K)$ is isomorphic to $\co$ for some set $\Gamma$;
\item[(iv)] $C(K)$ is isomorphic to a subspace  of
$\co$ for some set $\Gamma$.
\end{itemize}
\end{thm}

\begin{cor}\label{cor_images}
The class $\ecf$ is preserved under continuous images.
\end{cor}

\begin{proof}
Let $K\in\ecf$ and let $\phe: K\to L$ be a continuous surjection. By Theorem \ref{old_charact}, the space $C(K)$ is is isomorphic to $\co$ for some set $\Gamma$. The space $C(L)$ is isometric to a subspace $\{f\circ\phe: f\in C(L)\}$ of $C(K)$, hence it is isomorphic to a subspace of $\co$. Again, by Theorem \ref{old_charact}, the space $L$ belongs to $\ecf$.
\end{proof}

\begin{thm}\label{k(n)}
For each $n\in\omega$ there exists $k(n)\in\omega$ such that any continuous image of a space $K\in \ecn$ belongs to $\ec_{k(n)}$.
\end{thm}

\begin{proof}
Fix $n\in\omega$ and assume towards a contradiction, that, for each $i\in\omega$, there exist a compact space $K_i\in\ecn$ and a continuous surjection $\phe_i: K_i\to L_i$ such that $L_i$ does not belong to $\ec_i$. Without loss of generality we can assume that $K_i\subseteq \sigma_n({\Gamma_i})$ for some set $\Gamma_i$, and the sets $\Gamma_i$ are pairwise disjoint and disjoint with $\omega$. Consider 
$\Gamma = \omega\cup \bigcup_{i\in\omega} \Gamma_i\,.$ Let $X = \bigoplus _{i\in\omega} K_i$ and $Y = \bigoplus _{i\in\omega} L_i$ be discrete unions of spaces $K_i$ and $L_i$, respectively. Clearly, the one point compactification $\alpha(X)$ can be mapped continuously onto $\alpha(Y)$. Indeed, it is enough to take the union of all maps $\phe_i$, and assign $\infty_Y$ to $\infty_X$. Let $\psi: \alpha(X) \to \sigma_{n+1}(\Gamma)$ be defined by
$$\psi(x) = \begin{cases} \chi_{A\cup\{i\}}& \text{if $x= \chi_A\in
K_i$}\,,i\in\omega\,,\\
 \chi_\emptyset& \text{if $x = \infty_X$}\,, 
\end{cases}$$
for $x\in \alpha(X)$. A routine verification shows that $\psi$ is an embedding. On the other hand the compactification $\alpha(Y)$ does not belong to any class $\ec_i$, since these classes are hereditary with respect to closed subspaces. This gives a contradiction with Corollary \ref{cor_images}.
\end{proof}

\begin{ex}\label{ex_cont_image}
There exists a continuous image of the space $\sigma_2(\omega_1)$ which does not belong to $\ec_2$.
\end{ex}

\begin{proof}
Obviously, it is enough to construct an appropriate  continuous image of the space $\sigma_2(\omega_1\times 2)$ instead of $\sigma_2(\omega_1)$. 

Let $\sim$ be the equivalence relation on $\sigma_2(\omega_1\times 2)$ given by $\chi_{\{(\alpha,0)\}} \sim \chi_{\{(\alpha,1)\}}$, for all $\alpha\in \omega_1$, and let $q: \sigma_2(\omega_1\times 2) \to \sigma_2(\omega_1\times 2)/_{\sim}$ be the corresponding quotient map. Denote the quotient space $\sigma_2(\omega_1\times 2)/_{\sim}$ by $L$. It is routine to verify that the space $L$ is Hausdorff. We will show that $L\notin \ec_2$.

Suppose the contrary, then, by Lemma \ref{minimal_X}, we can assume that there exists an embedding $\phe: L\to \sigma_2(\omega_1)$. Since $L^{(2)} = \{q(\chi_{\emptyset})\}$ and $\sigma_2(\omega_1)^{(2)} = \{\chi_{\emptyset}\}$, we have $\phe(q(\chi_{\emptyset})) = \chi_{\emptyset}$. Therefore, the set $L^{(1)}\setminus L^{(2)}$ must be mapped by $\phe$ into $\sigma_2(\omega_1)^{(1)}\setminus \sigma_2(\omega_1)^{(2)}$. This means that there is an injective map $\psi: \omega_1\to \omega_1$ such that
\begin{equation}
\phe(q(\chi_{\{(\alpha,i)\}})) = \chi_{\{\psi(\alpha)\}}\quad \text{for all $\alpha\in \omega_1$ and $i=0,1$.}
\end{equation}
By the continuity of $\phe$ at the points $q(\chi_{\{(\alpha,i)\}})$ it follows that, for each $\alpha\in \omega_1$ there exits a finite set $F(\alpha)\subseteq \omega_1$ such that
\begin{equation}\label{ex_eq2}
(\forall \beta \in (\omega_1\setminus F(\alpha)))\  (\forall i,j\in 2)\ [\phe(\chi_{\{(\alpha,i), (\beta,j)\}}) \in V_{\{\psi(\alpha)\}}]
\end{equation}
(recall that $V_{\{\psi(\alpha)\}}$ is the clopen neighborhood $\{\chi_A\in \sigma_2(\omega_1): \psi(\alpha)\in A\}$ of $\chi_{\{\psi(\alpha)\}}$).
Take any $\gamma\in \omega_1\setminus \omega$ greater than  $\sup\bigcup_{n\in\omega} F(n)$. Observe that, for all $n\in\omega$, $\gamma\notin F(n)$. Pick any $k\in \omega\setminus F(\gamma)$. By property (\ref{ex_eq2}) we have
\begin{equation}
(\forall i,j\in 2)\ [\phe(\chi_{\{(k,i), (\gamma,j)\}}) \in V_{\{\psi(k)\}}\cap V_{\{\psi(\gamma)\}}]\,.
\end{equation}
Hence, the intersection $V_{\{\psi(k)\}}\cap V_{\{\psi(\gamma)\}}$ contains at least $4$ points (note that $k\ne \gamma$). On the other hand, for all distinct $\alpha,\beta\in \omega_1$, we have $V_{\{\alpha\}}\cap V_{\{\beta\}} = \{\chi_{\{\alpha,\beta\}}\}$, a contradiction.
\end{proof}

\begin{thm}\label{k(2)}
Each continuous image of a space $K\in \ec_2$ belongs to $\ec_{3}$.
\end{thm}

\begin{proof}
Let $\Gamma$ be a set, $K$ be a closed subset of $\sigma_2(\Gamma)$, and let $\phe: K\to L$ be a continuous surjection. We will show that $L$ embeds into $\sigma_3(\Gamma)$. Without loss of generality we can assume that the map $\phe$ is irreducible.

We will consider three cases determined by the height $ht(L)$ of $L$.\smallskip

\noindent
{\bf Case 1.} $ht(L)\le 1$. This means that $L$ is finite and this case is trivial.\smallskip

\noindent
{\bf Case 2.} $ht(L)=2$. In this case, for some $k\in\omega$, the space $L$ is homeomorphic to the discrete union $\bigoplus _{i=1}^k \alpha(X_i)$ of one point compactifications of infinite discrete spaces $X_i$, $i = 1,\dots,k$. Since $w(L)\le w(K)\le |\Gamma|$, we have $|X_i|\le |\Gamma|$ for $i = 1,\dots,k$.  Then $L$ embeds in $\sigma_2(\Gamma)$ by Lemma \ref{emb_n_into_n+1} and a simple observation that the space $\alpha(X_i)$ is homeomorphic to $\sigma_1(X_i)$.
\smallskip

\noindent
{\bf Case 3.} $ht(L)=3$. Then also $ht(K)=3$ and $K^{(2)} = \{\chi_\emptyset\} = \sigma_0(\Gamma)$. By Proposition \ref{scatt_Fact1}, $L^{(2)}$ is the singleton $\{\phe(\chi_\emptyset)\}$.  Since $\phe$ is irreducible, from Propositions \ref{scatt_Fact1} and  \ref{scatt_Fact2} it follows that
\begin{enumerate}[(a)]

\item $\phe\upharpoonright (K\setminus K')$  is a bijection onto $L\setminus L'$;

\item  $L'\setminus L^{(2)}\subseteq \phe(K'\setminus K^{(2)}) \subseteq L'$. 
\end{enumerate}

For each $y\in L'\setminus L^{(2)}$, the fiber $\phe^{-1}(y)$ is a closed in $K$ subset of $K'\setminus K^{(2)}$, hence it is finite. Since $K'\subseteq \sigma_1(\Gamma)$,  we have $\phe^{-1}(y)\subseteq \sigma_1(\Gamma)\setminus \sigma_0(\Gamma)$.

Recall that, for $\gamma\in\Gamma$, $V_{\{\gamma\}}$  denotes the clopen neighborhood $\{\chi_A\in \sigma_2(\Gamma): \gamma\in A\}$ of $\chi_{\{\gamma\}}$. We put  $U_{\{\gamma\}} = V_{\{\gamma\}}\cap K$. For $y\in L'\setminus L^{(2)}$ define
$$W_y = \bigcup \{U_{\{\gamma\}}: \chi_{\{\gamma\}}\in \phe^{-1}(y)\}\quad \text{and}\quad O_y = \phe\left(W_y\right).$$
Using properties (a) and (b) of $\phe$ one can easily verify that
\begin{enumerate}[(c)]
\item $O_y\cap L' = \{y\}$;
\end{enumerate}
\begin{enumerate}[(d)]
\item $\phe^{-1}(O_y) = W_y$ . 
\end{enumerate}
Clearly the set $W_y$ is clopen in $K$, since $\phe^{-1}(y)$ is finite. Therefore, by (d), $O_y$ is clopen in $L$. Let
$$\mathcal{U} = \{O_y: y\in L'\setminus L^{(2)}\}\cup \{\{z\}: z\in L\setminus L'\}\,.$$ 
The family $\mathcal{U}$ consists of clopen sets and by property (c) is $T_0$-separating in $L$. Let us check that this family is point-3. Since, for distinct $z_1, z_2\in L\setminus L'$, the singletons  $\{z_1\},\{z_2\}$ are obviously disjoint, it is enough to verify that, for distinct $y_1,y_2,y_3 \in L'\setminus L^{(2)}$, the intersection $\bigcap_{i=1}^3 O_{y_i}$ is empty. This follows from property (d) and an observation that for any $\chi_{\{\gamma_i\}}\in \phe^{-1}(y_i),\ i=1,2,3$ we have $\bigcap_{i=1}^3 V_{\{\gamma_i\}} = \emptyset$. Finally, our thesis follows from Proposition \ref{charact_ecn}.
\end{proof}

The following two easy observations demonstrate that the Example \ref{ex_cont_image} is in some sense the simplest possible.

\begin{prp}\label{metrizable_image}
For any $n\in\omega$, each continuous image of a metrizable space $K\in \ecn$ belongs to $\ecn$.
\end{prp}

\begin{proof}
By the classical characterization of countable compacta due to Mazur\-kie\-wicz and Sierpi\'nski \cite{MS}, any  countable metrizable space of height $n+1,\ n\in \omega$, is homeomorphic to the discrete union of $k$ copies of the ordinal space $\omega^n+1$, for some positive integer $k$. By Proposition \ref{scatt_Fact1}, any continuous image of such space is  homeomorphic to the discrete union of $k'$ copies of $\omega^{n'}+1$, where either $n'<n$ or $n'=n$ and $k'\le k$. Lemma \ref{minimal_X} implies that any metrizable space $K\in \ecn$ embeds into $\sigma_n(\omega)$. One can easily verify that, for any $n\in\omega$, the space $\sigma_n(\omega)$ is homeomorphic to the space $\omega^n+1$. The desired conclusion follows easily either from standard properties of ordinal spaces or Lemma \ref{emb_n_into_n+1}.
\end{proof}

The next proposition is trivial since the class $\ec_1$ consists of compact spaces with at most one nonisolated point.

\begin{prp}\label{k(1)}
Each continuous image of a space $K\in \ec_1$ belongs to $\ec_1$.
\end{prp}

Theorem \ref{k(2)} together with Example \ref{ex_cont_image} can be stated shortly that $3$ is the optimal value of the integer $k(n)$ from Theorem \ref{k(n)} for $n=2$. We do not know how to generalize this for $n> 2$.

\begin{prob}
Find the formula for best possible value of  $k(n)$ from Theorem \ref{k(n)}\footnote{Recently, Grzegorz Plebanek has proved the following recursive estimate for the optimal constant $k(n)$: $k(n)\le (2^n - 1)k(n-1) + 1$}.
\end{prob}

\begin{rem}\label{boolean_trans}
It is clear that the characterization of the class $\ecn$ from Proposition \ref{charact_ecn} can be formulated in the following way:

A compact zero-dimensional space $K$ belongs to $\ecn$, if and only if, the algebra $CO(K)$ of clopen subsets of $K$, has a set of generators $\mathcal{G}$ such that for any distinct $U_0,U_1,\dots,U_n\in \mathcal{G}$, the intersection $\bigcap_{i=o}^n U_i$ is empty. To simplify the statements, say for a moment, that the Boolean algebra $\mathcal{A}$ with such set of generators has the property $\mathcal{DG}_n$. Then Example \ref{ex_cont_image} can be translated into this language as follows: There exists a Boolean algebra $\mathcal{A}$ with property $\mathcal{DG}_2$, and a subalgebra $\mathcal{B}\subseteq \mathcal{A}$ without this property. Theorem \ref{k(2)} says that any subalgebra $\mathcal{B}$ of an algebra $\mathcal{A}$ with $\mathcal{DG}_2$, has the property $\mathcal{DG}_3$.
Other results from this section can be also reformulated in this way.
\end{rem}

\section{On zero-dimensional closed subspaces of nonmetrizable Eberlein compacta}

There are known several constructions, using some additional set-theoretic assumptions, of nonmetrizable compact spaces without nonmetrizable zero-dimen\-sional closed subspaces. Recently,
P.\ Koszmider \cite{Ko} constructed the first such example without such additional assumptions, and G. Plebanek \cite{Pl} constructed a consistent example of such a space which is a Corson compact. 
So it is important to determine whether we can obtain such examples within some other known classes of compact spaces, for example Eberlein compact spaces.
\smallskip

Joel Alberto Aguilar has asked us the following question. 

\begin{prob}\label{prob_refl_k}
Let $K$ be an Eberlein compact space of weight $\kappa$. Does $K$ contain a closed zero-dimensional subspace $L$ of the same weight?
\end{prob}

Probably, the most natural and interesting is the following simplified version of this question.

\begin{prob}\label{prob_refl}
Let $K$ be a nonmetrizable Eberlein compact space. Does $K$ contain a closed nonmetrizable zero-dimensional subspace $L$?
\end{prob}

We will show that the negative answer to this problem is consistent with \textsf{ZFC}, see Corollary \ref{ex_no_zerodim_sub}. We do not know if the affirmative answer is also consistent with \textsf{ZFC}, see Problem \ref{prob_refl_cons}.

We will also consider a more complex version of Problem \ref{prob_refl_k}:

\begin{prob}\label{prob_refl_k_l} Let  $\kappa \le \lambda$ be uncountable cardinal numbers, and
let $K$ be an Eberlein compact space of weight $\lambda$. Does $K$ contain a closed zero-dimensional subspace $L$ of weight $\kappa$?
\end{prob}

We will prove that, for every natural numbers  $1\le k\le n,\ n\ge 2$, the statement: each Eberlein compact space of weight $\omega_n$ contains a closed zero-dimensional subspace $L$ of weight $\omega_k$ is independent of \textsf{ZFC}, see Corollary \ref{ex_no_zerodim_sub_om_n} and Corollary \ref{cor_corson}.
\medskip

We begin with a simple observation based on the following well known property of Corson compacta. Since we were not able to find a reference for it, we include a sketch of a proof of this fact.

\begin{prp}\label{remark_character}
Let $x$ be a nonisolated point of a Corson compact space $K$ such that the character $\chi(K,x) = \kappa$. Then $K$ contains a copy of a one point compactification $\alpha(\kappa)$ of a discrete space of cardinality $\kappa$ with $x$ as its point
at infinity.
\end{prp}

\begin{proof}
Without loss of generality we can assume that $K\subseteq \Sigma(\Gamma)$, for some set $\Gamma$, and, for each $\gamma\in \Gamma$, there is $x_\gamma\in K$ such that $x_\gamma(\gamma) \ne 0$. Translating $K$ by the vector $-x$, we can also assume that $x= 0$ - the constant zero function in $\Sigma(\Gamma)$. If $\kappa = \omega$, then, by the Fr\'echet-Urysohn property of $K$, $0$ is the limit of a sequence of distinct points of $K$, which gives the desired conclusion. If  $\kappa > \omega$, then, using transfinite induction of length $\kappa$, and the fact that for compact spaces the pseudocharcter of a point is equal to the character, one can easily construct a set $\{x_\alpha: \alpha \in \kappa\}\subset K$ of points with nonempty, pairwise disjoint supports. Then the subspace $\{0\}\cup \{x_\alpha: \alpha \in \kappa\}$ of $K$ is as desired.
\end{proof}

\begin{cor}\label{cor_character_1}
Let $K$ be an Eberlein compact space with a point of character $\kappa$. Then $K$ contains a closed zero-dimensional subspace $L$ of weight $\kappa$. In particular, each Eberlein compact space of uncountable character contains a closed nonmetrizable zero-dimensional subspace $L$.
\end{cor}

It is worth to recall here that, by Arhangel'skii's theorem (\cite[Theorem 3.1.29]{En1}, for an infinite compact space $K$, we have the estimate $|K|\le 2^{\chi(K)}$. Hence we obtain the following.
\medskip

\begin{cor}\label{cor_character_2}
Let $K$ be an Eberlein compact space of weight greater that $2^\kappa$. Then $K$ contains a closed zero-dimensional subspace $L$ of weight $\kappa^+$. In particular, each Eberlein compact space of cardinality greater than continuum contains a closed nonmetrizable zero-dimensional subspace $L$.
\end{cor}

Recall that the definition of a $(\lambda,\kappa)$-Luzin set appearing in the next lemma, can be found in Subsection \ref{subs_Luzin}.

\begin{ex}\label{ex_no_zerodim_sub_k_l}  Let  $\kappa \le \lambda$ be uncountable cardinal numbers, and
assume that there exists a $(\lambda,\kappa)$-Luzin set. Then, for each $n\in\omega$ ($n=\infty$), there exists an $n$-dimensional nonmetrizable Eberlein compact space $K_n$ of weight $\lambda$ such that any closed subspace $L$ of $K_n$ of weight $\ge \kappa$  has dimension $n$.
\end{ex}

Applying the above for the standard Luzin set we obtain that it is consistent that Problem \ref{prob_refl} (Problem \ref{prob_refl_k}) has an negative answer.

\begin{cor}\label{ex_no_zerodim_sub}
Assuming the existence of a Luzin set, there exists a nonmetrizable Eberlein compact space $K$ without closed nonmetrizable zero-dimensional subspaces.
\end{cor}

As we mentioned in Subsection \ref{subs_Luzin}, for each $n\ge 1$, the existence of a $(\omega_n,\omega_1)$-Luzin set in $\mathbb{R}$ is consistent with \textsf{ZFC}, therefore we derive the following.

\begin{cor}\label{ex_no_zerodim_sub_om_n} For each $n\ge 1$, it is consistent with \textsf{ZFC} that there exists an Eberlein compact space $K$ of weight $\omega_n$ without closed nonmetrizable zero-dimensional subspaces.
\end{cor}

\begin{proof}[Construction of Example \ref{ex_no_zerodim_sub_k_l}]
Fix $n\in\omega$ ($n=\infty$). Let $X$ be a $(\lambda,\kappa)$-Luzin set in the cube $[0,1]^n$, see Subsection \ref{subs_Luzin}. We consider the following subspace of the Aleksandrov duplicate $AD([0,1]^n)$ (see Subsection \ref{subs_AD})
$$K = AD([0,1]^n)\setminus [([0,1]^n\setminus X)\times \{1\}]\,.$$
Since $\dim(AD([0,1]^n)) = n$ and $K$ contains a copy of the cube $[0,1]^n$, the compact space $K$ is $n$-dimensional (cf.\ \cite[Ch.\ 3]{En}). By Proposition \ref{prp_AD_Eberlein}, $K$ is uniform Eberlein compact. 

Let $L$ be a closed subspace of $K$ of  weight $\ge \kappa$ and let
$$Z = L\cap ([0,1]^n\times \{1\}) = L\cap ([X\times \{1\}) = Y\times \{1\}\,.$$
Since $L$ has  weight $\ge \kappa$, the set $Y$ must have the cardinality $\ge \kappa$. Let $T\subseteq Y$ be a subset of the same cardinality as $Y$, and without isolated points. By the definition of a $(\lambda,\kappa)$-Luzin set, $T$ is not nowhere dense in $[0,1]^n$, hence its closure $\clo_{[0,1]^n}(T)$ has a nonempty interior in $[0,1]^n$. Therefore, $\clo_{[0,1]^n}(T)$ has dimension $n$. Since $T$ is dense-in-itself, the closure of $T\times \{1\}$ in $AD([0,1]^n)$ (hence also in $L$) contains the set $\clo_{[0,1]^n}(T)\times \{0\}$, a topological copy of $\clo_{[0,1]^n}(T)$. It follows that $\dim (L) = n$.
\end{proof}

Now, we will switch to the consistent results giving an affirmative answer to some cases of Problem \ref{prob_refl_k_l}.

We start with the following technical lemma (the definitions of cardinal numbers used in this lemma can be found in Subsection \ref{subs_b_nonM}).

\begin{lem}\label{lem_almost_disj_supp}
Let $\kappa$ be a cardinal number of uncountable cofinality and assume that $\mathfrak{b} > \kappa$. Let $K$ be a compact subset of the product $\mathbb{R}^\Gamma$ containing a subset $X$ of cardinality $\kappa$ such that, for some countable subset $\Gamma_0$ of $\Gamma$ and for all $x\in X$, the sets $\supp(x)\setminus \Gamma_0$ are nonempty and pairwise disjoint. Then $K$ contains a closed zero-dimensional subspace $L$ of weight $\kappa$.
\end{lem}

\begin{proof}
From our assumption that  $\mathfrak{b} > \kappa$, it follows that also $\mathrm{non}(\mathcal{M}) > \kappa$.

Let $\Gamma_0\subseteq \Gamma$ be a countable set witnessing the property of the set $X$. Observe that the property of supports of points $x\in X$, implies that $X$ is a discrete subspace of $K$.

Without loss of generality we can assume that the set $\Gamma_0$ is infinite, so we can enumerate it as $\{\gamma_n: n\in\omega\}$.
Let $X_n = \{x(\gamma_n): x\in X\}$ for $n\in\omega$. Since $\mathrm{non}(\mathcal{M}) > \kappa$, each set $X_n$ is meager. Therefore, for each $n\in\omega$, we can find an increasing sequence $(C_n^k)_{k\in\omega}$ of closed nowhere dense subsets of $\mathbb{R}$ such that $X_n\subseteq \bigcup_{k\in\omega} C_n^k$. For each $x\in X$, we define a function $f_x: \omega\to \omega$ as follows
$$f_x(n) = \min\{k: x(\gamma_n)\in C_n^k\}\quad \text{for } n\in\omega.$$
Since $\mathfrak{b} > \kappa$, we can find a function $g: \omega\to \omega$ such that $f_x\le^* g$ for all $x\in X$.  A routine refining argument, using uncountable cofinality of $\kappa$, shows that there is a subset $Y\subseteq X$ of size $\kappa$ and a function $h: \omega\to \omega$ such that $f_x\le h$ for all $x\in Y$. 

We define $L = \clo_K Y$.

The space $L$ contains a discrete subspace $Y$ of cardinality $\kappa$, hence $L$ has weight $\ge \kappa$. 
Let $\Gamma_1$ be a selector from the family $\{\supp(x)\setminus \Gamma_0: x\in Y\}$. Clearly, $\Gamma_1$ has cardinality $\kappa$. One can easily verify that the projection $p: \mathbb{R}^\Gamma\to \mathbb{R}^{\Gamma_0\cup \Gamma_1}$ is one-to-one on $L$, hence $w(L)\le \kappa$.

It remains to verify that $L$ is zero-dimensional.

First, observe that each space $C_n^k$, being closed nowhere dense in $\mathbb{R}$ is zero-dimensional. Let $\pi: L\to \mathbb{R}^\omega$ be defined by
$$\pi(x)(n) = x(\gamma_n)\quad \text{for } x\in L,\  n\in\omega.$$
From our choice of $Y$ and $h$ it follows that $$\pi(L)\subseteq \Pi_{n\in\omega} C_n^{h(n)}.$$
Since the product $\Pi_{n\in\omega} C_n^{h(n)}$ is zero-dimensional, so is the space $\pi(L)$. From the fact that the sets $\supp(x)\setminus \Gamma_0$, for $x\in Y$, are pairwise disjoint, it follows that each fiber of $\pi$ is either finite, or homeomorphic to a one point compactification of a discrete space, hence it is zero-dimensional. Therefore, by the theorem on dimension-lowering mappings \cite[Theorem 3.3.10]{En}, $L$ is zero-dimensional.
\end{proof}

\begin{lem}\label{lem_count_1}
Let $k<n$ be natural numbers, $Y$ be a set of cardinality $\omega_k$, and $\mathcal{C} = \{C_\alpha: \alpha < \omega_n\}$ be a family of countable subsets of $Y$. Then there exist a countable subset $Z$ of $Y$ and a subset $S$ of $\omega_n$ of cardinality $\omega_n$ such that $C_\alpha\subseteq Z$ for all $\alpha\in S$.
\end{lem}

\begin{proof}
Fix $n\ge 1$. The case $k=0$ is trivial. For $k>0$, we proceed by induction on $k$. Without loss of generality we can assume that $Y= \omega_k$. Since the cofinality of $\omega_k$ is uncountable, we can find an $\lambda<\omega_k$ and a subset $R$ of $\omega_n$ of cardinality $\omega_n$ such that $C_\alpha \subseteq \lambda$ for $\alpha \in R$. Now, we can use the inductive hypothesis.
\end{proof}

The following lemma is probably well known. We learned about it from Grzegorz Plebanek, who suggested to use it for the proof of Theorem \ref{thm_corson}. Its proof is based on an idea from the proof of theorem 1.6 in \cite{Ku}.

\begin{lem}\label{lem_count_2}
Let $\Gamma$ be a set of cardinality $\omega_n$, $n\ge 2$, and $\mathcal{A} = \{A_\alpha: \alpha < \omega_n\}$ be a family of countable subsets of $\Gamma$, such that $\bigcup \mathcal{A} = \Gamma$. Then there exist a countable subset $B$ of $\Gamma$ and a subset $T$ of $\omega_n$ of cardinality $\omega_n$ such that the family  $\{A_\alpha\setminus B: \alpha \in T\}$ consists of nonempty, pairwise disjoint sets.
\end{lem}

\begin{proof}
Without loss of generality we can assume that $\Gamma = \omega_n$. For each $\beta  < \omega_n$ we pick an $\alpha(\beta) < \omega_n$ such that $\beta \in A_{\alpha(\beta)}$.
We consider the sets $A_\alpha$ with the order inherited from $\omega_n$. Since we have only $\omega_1$ possible order types of these sets, we can find a countable ordinal $\eta$, and a subset $P$ of $\omega_n$ of cardinality $\omega_n$ such that, for all $\beta \in P$, $A_{\alpha(\beta)}$ has order type $\eta$. For $\gamma < \eta$, let $\xi(\beta,\gamma)$ be the $\gamma$-th element of $A_{\alpha(\beta)}$. Since the union $\bigcup\{A_{\alpha(\beta)}: \beta \in P\}$ contains an unbounded set $P$, and $\omega_n$ has uncountable cofinality, there is $\gamma < \eta$ such that the set $\{\xi(\beta,\gamma): \beta \in P\}$ is unbounded in $\omega_n$. Let $\gamma_0$ be the smallest such $\gamma$. Put
$$\delta = \sup\{\xi(\beta,\gamma): \beta \in P,\ \gamma<\gamma_0\} + 1\,.$$
Using the definition of $\gamma_0$, one can easily construct, by a transfinite induction of length $\omega_n$, a subset $Q$ of $P$ of cardinality $\omega_n$, such that the family  $\{A_{\alpha(\beta)}\setminus \delta: \beta \in Q\}$ consists of nonempty, pairwise disjoint sets. Finally, we can apply Lemma \ref{lem_count_1}, for $Y= \delta$ and the family $\{A_{\alpha(\beta)}\cap \delta: \beta \in Q\}$, to find a countable subset $B$ of $\delta$ and a subset $T$ of $Q$ of cardinality $\omega_n$, such that $(A_{\alpha(\beta)}\cap \delta)\subseteq B$ for $\beta\in T$. 
\end{proof}

\begin{thm}\label{thm_corson}
Assume that $\mathfrak{b} > \omega_n,\ n\ge 1$. Then each Corson compact space $K$ of weight greater that $\omega_1$ contains a closed zero-dimensional subspace $L$ of weight equal to $\min(w(K), \omega_n)$.
\end{thm}

\begin{cor}\label{cor_corson} For each $n\ge 2$, it is consistent with \textsf{ZFC} that each Corson compact space $K$ of weight $\omega_n$ contains a closed zero-dimensional subspace $L$ of the same weight.
\end{cor}

\begin{proof}[Proof of Theorem \ref{thm_corson}]
Let $\lambda = w(K)\ge \omega_2$ and $\kappa = \min(w(K), \omega_n)$. Let $\eta = \omega_2$ if $\kappa = \omega_1$, otherwise $\eta = \kappa$. Without loss of generality we can assume that $K\subseteq \Sigma(\lambda)$ and, for each $\gamma\in \lambda$, there is $x_\gamma\in K$ such that $x_\gamma(\gamma) \ne 0$. Pick a subset $S$ of $\lambda$ of cardinality $\eta$ and put $\Gamma = \bigcup\{\supp(x_\gamma): \gamma \in S\}$. We apply Lemma \ref{lem_count_2} for $\Gamma$ and the family $\{\supp(x_\gamma): \gamma \in S\}$, to obtain a countable subset $B$ of $\Gamma$ and subset $T\subseteq S$ of cardinality $\eta$ such that the family  $\{\supp(x_\gamma)\setminus B: \gamma \in T\}$ consists of nonempty, pairwise disjoint sets. If $\kappa > \omega_1$, take $X = \{x_\gamma: \gamma\in T\}$, otherwise pick a subset $T_0\subseteq T$ of cardinality $\omega_1$ and put $X = \{x_\gamma: \gamma\in T_0\}$. Now, we can obtain the desired conclusion applying Lemma \ref{lem_almost_disj_supp}.
\end{proof}

\medskip

The construction from Example \ref{ex_no_zerodim_sub_k_l} and Lemma \ref{lem_almost_disj_supp} motivated us to consider the following class of Eberlein compacta.

We say that a compact space $K$ belongs to the class $\eco$ if, for some set $\Gamma$ there is an embedding $\phe:K\to \mathbb{R}^\Gamma$ and a countable subset $\Gamma_0$ of $\Gamma$ such that, for each $x\in K$, the set $\supp(\phe(x))\setminus \Gamma_0$ is finite. Since the product $\mathbb{R}^{\Gamma_0}$ embeds into the Hilbert space $\ell_2$ equipped with the pointwise topology, it easily follows that any compact space $K\in\eco$ is  uniform Eberlein compact. It is clear that the class $\eco$ is preserved by the operations of taking finite products and closed subspaces. Example \ref{ex_non_eco} below demonstrates that the countable power of a space from $\eco$ may not belong to this class.

One can also easily verify that, for a metrizable compact space $M$ the Aleksandrov duplicate $AD(M)$ belongs to the class $\eco$ (cf.\ the proof of Proposition \ref{prp_AD_Eberlein}), hence all spaces constructed in such a way as in Example \ref{ex_no_zerodim_sub_k_l}, are in this class. 
One can even show that, for metrizable compacta $M_n$, the countable product of spaces $AD(M_n)$ is in $\eco$. In particular, the product $[AD(2^\omega)]^\omega$ belongs to $\eco$. Note, that by the remarkable result of Dow and Pearl \cite{DP} this product is an example of a homogeneous nonmetrizable Eberlein compact space. The first such example was given by Jan van Mill in \cite{vM}. The structure of these two examples seems to be closely related, but we do not know if they are homeomorphic.

\begin{thm}\label{thm_on_eco}
Assuming that $\mathfrak{b} > \omega_1$, each nonmetrizable compact space $K\in\eco$ contains a closed nonmetrizable zero-dimensional subspace $L$.
\end{thm}

\begin{proof}
Let $K\in\eco$. Without loss of generality we can assume that, for some set $\Gamma$ and its countable subset $\Gamma_0$, $K$ is a subset of  $\mathbb{R}^\Gamma$ such that, for each $x\in K$, the set $\supp(x)\setminus \Gamma_0$ is finite. Since $K$ is nonmetrizable, obviously the set $\Gamma$ must be uncountable. We can also assume that, for each $\gamma\in \Gamma$, there is $x_\gamma\in K$ such that $x_\gamma(\gamma) \ne 0$. For each $\gamma\in \Gamma\setminus \Gamma_0$, the set $F_\gamma = \supp(x_\gamma)\setminus \Gamma_0$ is finite and nonempty. Using the $\Delta$-system lemma we can find a finite set $A\subseteq \Gamma$ and a set $S\subseteq (\Gamma\setminus \Gamma_0)$ of size $\omega_1$ such that, for any distinct $\alpha,\beta \in S$, $F_\alpha\cap F_\beta = A$. By enlarging $\Gamma_0$ to $\Gamma_0\cup A$, we can assume that $A$ is empty. Now, we can apply Lemma \ref{lem_almost_disj_supp} for the set $X=\{x_\gamma: \gamma\in S\}$. 
\end{proof}

In the light of Corollary \ref{cor_character_1} and Theorem \ref{thm_on_eco} it seems natural to ask whether every first-countable Eberlein compact space belongs to the class $\eco$. Unfortunately, this is not the case.

\begin{ex}\label{ex_non_eco}
There exists a first-countable uniform Eberlein compact space which does not belong to the class $\eco$.
\end{ex}

\begin{proof}
Our construction uses the following modification $L$ of the Aleksandrov duplicate $AD([0,1])$ of the unit interval (roughly speaking, we replace isolated points of $AD([0,1])$ by copies of $[0,1]$). The space $L$ is similar to the space $X(C)$ used by van Mill in \cite{vM} to construct an example of a homogeneous nonmetrizable Eberlein compact space.

Let $S= [0,1]\cup\{2\}$. For any $t\in [0,1],\ u\in [1,2]$, we define functions $f_t: S\to [0,2],\ g_{t,u}: S\to [0,2]$ by
\begin{eqnarray*}
f_t(s) &=& \begin{cases} 0\quad \text{if } s\in [0,1]\,,\\
t\quad \text{if } s = 2\,; 
\end{cases}\\
g_{t,u}(s) &=& \begin{cases} u\quad \text{if } s=t\,,\\
0\quad \text{if } s\in [0,1],\ s\ne t\,,\\
t\quad \text{if } s = 2\,. 
\end{cases}\\
\end{eqnarray*}
We consider 
$$L = \{f_t: t\in[0,1]\}\cup \{g_{t,u}: t\in [0,1],\ u\in [1,2]\}$$
as a subspace of the cube $[0,2]^S$. One can easily verify that $L$ is closed in $[0,2]^S$. Since the cardinality of supports of functions $f_t$ and $g_{t,u}$ is bounded by $2$, $L$ is an uniform Eberlein compact space. It is also easy to observe that the space $L$ is first-countable. 

We will show that the space $K= L^\omega$ has the required property. Clearly, it is enough to show that $K$ does not belong to the class $\eco$.

Suppose the contrary, i.e., there exist a set $\Gamma$, its countable subset $\Gamma_0$, and an embedding $\phe:K\to \mathbb{R}^\Gamma$  such that, for each $x\in K$, the set $\supp(\phe(x))\setminus \Gamma_0$ is finite. We will treat $K$ as a subset of the product $[0,2]^{\omega\times S}$, namely we identify the sequence $(x_n)_{n\in\omega}\in  L^\omega$ with the function $x: \omega\times S$ defined by $x(n,s) = x_n(s)$ for $n\in\omega$, $s\in S$. Let $\pi: \mathbb{R}^\Gamma\to \mathbb{R}^{\Gamma_0}$ be the projection, and let $\psi = \pi\circ \phe$. By the Tietze extension theorem we can extend $\psi$ to a continuous map $\Psi: [0,2]^{\omega\times S}\to \mathbb{R}^{\Gamma_0}$. It is well known that such a map depends on countably many coordinates, i.e., there is a countable subset $T\subseteq \omega\times S$ and a continuous map $\theta: [0,2]^T\to \mathbb{R}^{\Gamma_0}$ such that $\Psi = \theta\circ \rho$, where $\rho$ denotes the projection of $[0,2]^{\omega\times S}$ onto $[0,2]^T$ (cf.\  \cite[2.7.12]{En1}). Denote the restriction $\rho\upharpoonright K$ by $\upsilon$. Clearly, $\theta\circ\upsilon = \psi = \pi\circ \phe$.
This implies that, for any $y\in \psi(K)$, we have $\upsilon^{-1}(\theta^{-1}(y)) = \phe^{-1}(\pi^{-1}(y))$. In particular, this means that the sets $\upsilon^{-1}(\theta^{-1}(y))$ and  $\pi^{-1}(y)\cap \phe(K)$ are homeomorphic.
Observe that the set $\pi^{-1}(y)\cap \phe(K)$ can be treated as a subspace of the space $\sigma(\mathbb{R}^{\Gamma\setminus \Gamma_0})$, consisting of functions from $\mathbb{R}^{\Gamma\setminus \Gamma_0}$ with finite supports. The space $\sigma(\mathbb{R}^{\Gamma\setminus \Gamma_0})$ is \emph{strongly countable-dimensional}, i.e, is a countable union of closed finite-dimensional subspaces, cf.\ \cite[proof of Proposition 1]{EP}. Therefore, the space $\pi^{-1}(y)\cap \phe(K)$ is also strongly countable-dimensional.

Pick a point $t\in [0,1]$ such that $(\omega\times \{t\})\cap T = \emptyset$. Let $x = (x_n)_{n\in\omega}\in K$ be the constant sequence, where $x_n= f_t$ for all $n$, and let $y = \psi(x)$. One can easily verify that the set $\upsilon^{-1}(\upsilon(x))\subseteq \upsilon^{-1}(\theta^{-1}(y))$, contains the product $\{g_{t,u}: u\in [1,2]\}^\omega$ homeomorphic to the Hilbert cube $[0,1]^\omega$. Since the Hilbert cube is not strongly countable-dimensional (which follows easily from the Baire Category Theorem), the subspace $\upsilon^{-1}(\theta^{-1}(y))$ is not strongly countable-dimensional, a contradiction.
\end{proof}

Let us note that using a very similar argument as above one can show that the unit ball $B$ of the Hilbert space $\ell_2(\omega_1)$ equipped with the weak topology is an example of an uniform Eberlein compact space which does not belong to the class $\eco$. This is much simpler example than above one, but it is not first-countable.
\medskip

Corollaries \ref{ex_no_zerodim_sub_om_n} and  \ref{cor_corson} indicate that probably the most interesting and left open case of Problem \ref{prob_refl_k_l} is the following

\begin{prob}\label{prob_refl_cons}
Is it consistent that every Eberlein compact space $K$ of weight $\omega_1$ contains a closed  zero-dimensional subspace $L$ of the same weight?
\end{prob}

\subsection*{Acknowledgments}
We would like to thank Grzegorz Plebanek for many fruitful discussions on the topic of this paper and several valuable suggestions.


\begin{thebibliography}{En2}

\bibitem[Av]{Av} A.\ Avil\'es, \emph{Countable products of spaces of finite sets}, Fund. Math. \textbf{186} (2005), 147--159.

\bibitem[BJ]{BJ} T.\ Bartoszyński and H.\ Judah, \emph{Set Theory. On the Structure of the Real Line},  A K Peters, Ltd., Wellesley, MA, (1995).

\bibitem[Be]{Be} M.\ Bell, {\em A Ramsey theorem for polyadic spaces}, Fund. Math.
\textbf{150} (1996), 189-195

\bibitem[BM]{BM}
M.\ Bell and W.\ Marciszewski, \emph{On scattered Eberlein compact
spaces}, Israel J. Math. \textbf{158} (2007), 217--224.

\bibitem[vD]{vD} E.\ van Douwen,  \emph{The integers and topology, Handbook of Set-Theoretic Topology} (K. Kunen
and J. Vaughan, eds.), North-Holland, 1984, pp.\ 111--167.

\bibitem[DP]{DP} A.\ Dow and E.\ Pearl, \emph{Homogeneity in powers of zero-dimensional first-countable spaces}, Proc. Amer. Math. Soc. \textbf{125} (1997), 2503--2510.

\bibitem[En1]{En1}
R.\ Engelking, \emph{General Topology}, Heldermann Verlag, Berlin (1989).


\bibitem[En2]{En}
R.\ Engelking, {\em Theory of Dimensions Finite and Infinite}, Sigma Series in Pure Mathematics, 10. Heldermann Verlag, Lemgo, 1995.

\bibitem[EP]{EP} R.\ Engelking and R.\ Pol {\em Compactifications of countable-dimensional and strongly countable-dimensional spaces},
Proc. Amer. Math. Soc. \textbf{104} (1988), 985--987.

\bibitem[Ko]{Ko} P.\ Koszmider, {\em  On the problem of compact totally disconnected reflection of nonmetrizability}, Topology Appl. \textbf{213} (2016), 154--166.

\bibitem[KM]{KM} M.\ Krupski and W. Marciszewski, \emph{Some remarks on universality properties of $l_{\infty}/c_0$}, Coll. Math. \textbf{128} (2) (2012), 187--195.

\bibitem[Ku]{Ku} K.\ Kunen, \emph{Set Theory.
An Introduction to Independence Proofs}, Studies in Logic and the Foundations of Mathematics, 102. North-Holland Publishing Co., Amsterdam-New York, 1980. 

\bibitem[Ma]{Ma} W.\ Marciszewski, \emph{On Banach spaces $C(K)$
isomorphic to $c_0(\Gamma)$}, Studia Math. \textbf{156} (2003),
295--302.

\bibitem[MS]{MS} S.\ Mazurkiewicz and W.\ Sierpi\'nski,  \emph{Contribution \`{a} la topologie des ensembles d\'enombrables}, Fund. Math. \textbf{1} (1920), 17--27.

\bibitem[vM]{vM} J.\ van Mill, \emph{A homogeneous Eberlein compact space which is not metrizable}, Pacific J. Math. \textbf{101} (1982) 141--146.

\bibitem[Ne]{Ne} S. Negrepontis, \emph{Banach spaces and topology, 
Handbook of Set-Theoretic Topology} (K. Kunen
and J. Vaughan, eds.), North-Holland, 1984, pp.\  1045--1142

\bibitem[Pl]{Pl} G.\ Plebanek, \emph{Musing on Kunen's compact $L$-space}, preprint, 	arXiv:2012.01849.

\bibitem[Se]{Se}
Z.\ Semadeni, \textit{Banach Spaces of Continuous Functions}, PWN, Warsaw, 1971.


\end{thebibliography}
\end{document}